\documentclass{amsart}
\usepackage{amssymb}
\usepackage{hyperref}
\usepackage[utf8]{inputenc}
\usepackage{dsfont}

\usepackage{graphicx}     
\usepackage{hyperref}
\usepackage{cleveref}

\hypersetup{
    pdftoolbar=true,
    pdfmenubar=true,
    pdffitwindow=false,
    pdfstartview={FitH},
    pdftitle={},
    pdfauthor={},
    pdfsubject={},
    pdfkeywords={}
    pdfnewwindow=true,
    colorlinks=true, 
    linkcolor=blue,
    citecolor=blue,
    urlcolor=black,
}

\usepackage{color}

\newtheorem{theorem}{Theorem}[]

\newtheorem{question}{Question}

\newtheorem*{mainthm*}{\sc Theorem}

\theoremstyle{definition}

\theoremstyle{remark}
\newtheorem{remark}[theorem]{Remark}

\allowdisplaybreaks
\numberwithin{equation}{section}

\usepackage{dsfont}

\newcommand{\Z}{\mathds Z}

\newcommand{\R}{\mathds R}
\newcommand{\C}{\mathds C}

\DeclareMathOperator{\Ot}{\mathsf{O}}
\DeclareMathOperator{\SO}{\mathsf{SO}}

\newcommand{\Id}{\textup{Id}}

\DeclareMathOperator{\diam}{diam}

\DeclareMathOperator{\dist}{d}

\newcommand{\G}{\mathsf{G}}

\newcommand{\Ut}{\mathsf{U}}

\newcommand{\Ss}{\mathds S}
\newcommand{\g}{\mathrm g}
\newcommand{\dd}{\mathrm d}

\allowdisplaybreaks
\numberwithin{equation}{section}

\title{Diameter and displacement of sphere involutions}

\author{Renato G. Bettiol}
\address{\!\!\!\begin{tabular}{lll}
CUNY Lehman College & & CUNY Graduate Center \\
Department of Mathematics & & Department of Mathematics \\
250 Bedford Park~Blvd W & & 365 Fifth Avenue \\
Bronx, NY, 10468, USA & & New York, NY, 10016, USA
\end{tabular}
}
\email{r.bettiol@lehman.cuny.edu}

\author{Emilio~A.~Lauret}
\address{Instituto de Matemática (INMABB)\newline\indent
Departamento de Matemática\newline\indent
Universidad Nacional del Sur (UNS)-CONICET\newline\indent
Bahía Blanca, Argentina}
\email{emilio.lauret@uns.edu.ar}

\date{\today}

\begin{document}

\begin{abstract}
We show that spheres in all dimensions $\geq3$ can be deformed to have diameter larger than the distance between any pair of antipodal points. This answers a question of Yurii Nikonorov. 
\end{abstract} 

\maketitle

\section{Introduction}

The diameter $\diam(M,\dist)$ of a compact length space is the maximal distance between pairs of points in $(M,\dist)$; if $M$ is a manifold and $\dist=\dist_\g$ is induced by a Riemannian metric $\g$, we write $\diam(M,\g)=\diam(M,\dist_\g)$. For example, the round $n$-sphere of radius $r$ has $\diam(\Ss^n(r))=\pi\,r$.
Nikonorov~\cite{Nikonorov01diameter} proved the following:

\begin{theorem}[Nikonorov]\label{thm:Yurii}
If $(M,\dist)$ is a length space homeomorphic to the sphere~$\Ss^2$ and $I\colon M\to M$ is an isometric involution without fixed points, then there exists $x\in M$ such that $\diam(M,\dist)=\dist(x,I(x))$. 
\end{theorem}

The above naturally leads to the following question \cite{Nikonorov01diameter}:

\begin{question}[Nikonorov]\label{question}
Is there an analogue of Theorem~\ref{thm:Yurii} for length spaces homeomorphic to the sphere $\Ss^n$ for some $n\geq 3$?
\end{question}

Podobryaev~\cite{Podobryaev18} observed that sufficiently collapsed Berger spheres provide a negative answer in dimension $n=3$. In fact, this observation can be easily extended to all \emph{odd} dimensions $n\geq3$, considering the (homogeneous) spheres $(\Ss^{2q+1},\g(t))$ obtained scaling the unit round sphere by $t>0$ in the vertical direction of the Hopf bundle $\Ss^1\to \Ss^{2q+1}\to \C P^q$. 
For all $t>0$, the projection onto $\C P^q$ remains a Riemannian submersion, hence $\diam(\Ss^{2q+1},\g(t)) \geq \diam(\C P^q) = \frac\pi2$. Meanwhile, pairs of antipodal points $x$ and $I(x)=-x$ on $(\Ss^{2q+1},\g(t))$ are also antipodal points on the totally geodesic fiber $\Ss^1(t)$, and thus $\dist_{\g(t)}(x,I(x))\leq \pi\, t$. Therefore, $\dist_{\g(t)}(x,I(x))<\diam(\Ss^{2q+1},\g(t))$ for all $t<\frac12$.
The latter actually holds for all $t<\frac{1}{\sqrt2}$ due to the explicit computation \eqref{eq:diam} of $\diam(\Ss^{2q+1},\g(t))$ by Rakotoniaina~\cite{Ra85}, recently rediscovered (in dimension 3) by Podobryaev~\cite{Podobryaev17}.

\medskip

In this short note, we provide similar negative answers in \emph{all dimensions $n\geq3$}.

\medskip

Our first construction involves the \emph{spherical join} $\Ss^k(r) * \Ss^{n-k-1}(r)$, $1\leq k \leq n-2$, of spheres of radius $0<r<\frac12$, which is a length space (in fact, an Alexandrov space) with diameter~$\frac\pi2$ and homeomorphic to~$\Ss^n$, see \cite[p.~582]{GrovePetersen} or \cite[p.~63]{bh} for definitions. Every point in $\Ss^k(r) * \Ss^{n-k-1}(r) \setminus (\Ss^k(r) \cup \Ss^{n-k-1}(r))$ has coordinates $(x,\rho,y)$, where $x\in \Ss^k(r)$, $y\in \Ss^{n-k-1}(r)$, and $\rho\in \left(0,\frac\pi2\right)$. There is a natural isometric action of $\SO(k+1)\times \SO(n-k)$ given by $(A,B)\cdot(x,\rho,y)=(Ax,\rho,By)$, whose orbits have diameter $\pi\, r<\frac\pi2$, since (see, e.g.,~\cite[p.~63]{bh}),
\begin{equation*}
     \dist^{\textrm{sph}}_{\textrm{join}}\!\big( (x_1,\rho,y_1),(x_2,\rho,y_2)\big) =\arccos\left(\cos^2\rho \,\cos(\dist(x_1,x_2))+\sin^2\rho \,\cos(\dist(y_1,y_2))\right),
\end{equation*}
which is bounded from above by $\max \{\dist(x_1,x_2), \dist(y_1,y_2) \}\leq \pi r$, where $\dist$ is used for distances in $\Ss^k(r)$ and $\Ss^{n-k-1}(r)$.
The involution $I(x,\rho,y)=(-x,\rho,-y)$ induced by the antipodal maps of each sphere is an isometry without fixed points, and corresponds to the antipodal map of $\Ss^n$ under the above homeomorphism.
Since $I$ commutes with the $\SO(k+1)\times \SO(n-k)$-action, it leaves invariant each orbit, and thus its maximal displacement is $ \pi\, r<\frac\pi2$. Therefore, $\Ss^k(r) * \Ss^{n-k-1}(r)$, with $1\leq k \leq n-2$ and $0<r<\frac12$, yields a negative answer to \Cref{question} for all $n\geq3$.

The spherical join $\Ss^k(r) * \Ss^{n-k-1}(r)$ is a smooth Riemannian manifold if and only if $r=1$, in which case it is isometric to $\Ss^n(1)$. However, inspired by this construction, we can also produce \emph{smooth} counter-examples to \Cref{question}, as follows:

\begin{mainthm*}
For all $n\geq3$, there is a family of smooth Riemannian metrics $(\g_s)_{s\geq 0}$ on~$\Ss^n$, such that $\g_0$ is the unit round metric, $\diam(\Ss^n,\g_s)\geq  \frac\pi2$, and the antipodal map $I(x)=-x$ is an isometry satisfying $\dist_{\g_s}(x,I(x))\leq\frac{\pi}{\sqrt{1+\frac{s}{2}}}$ for all~$x\in\Ss^n$.
\end{mainthm*}

Clearly, for $s>6$, the spheres $(\Ss^n,\g_s)$ provide a negative answer to \Cref{question} in all dimensions $n\geq3$.
These spheres are Cheeger deformations of $\Ss^n(1)\subset\R^{n+1}$ with respect to the block diagonal subgroup of isometries $\SO(k+1)\times \SO(n-k)$ in $\SO(n+1)$, with $1\leq k\leq n-2$. In particular, they are cohomogeneity one manifolds with geometric features similar to $\Ss^k(r) * \Ss^{n-k-1}(r)$; e.g., both are positively curved and converge in Gromov--Hausdorff sense to $\left[0,\frac\pi2\right]$ as $s\nearrow+\infty$, respectively $r\searrow0$. 
In fact, the unifying feature of all constructions in this note is that they are spheres with a distance-nonincreasing map onto $\left[0,\frac\pi2\right]$ whose fibers are invariant under the antipodal map and can be deformed to have arbitrarily small intrinsic diameter.

\subsection{Acknowledgements}
This paper is a contribution to the special issue ``TYAN Virtual Thematic Workshop in Mathematics'' of \textsl{Matem\'atica Contempor\^anea}. We thank the organizers for the invitation and their excellent job.
We would also like to thank Alberto Rodr\'iguez-V\'azquez for discussions about Berger spheres and for informing us of the paper \cite{Ra85}, Yurii Nikonorov for bringing Question~\ref{question} to our attention and for useful comments on a first draft of the paper, and the anonymous referee for thoughtful suggestions to improve the presentation. 

The first-named author is supported by the NSF CAREER grant DMS-2142575 and NSF grant DMS-1904342. The second-named author is supported by grants from FONCyT (PICT-2018-02073 and PICT-2019-01054) and SGCYT--UNS.

\section{Main construction}

Let $\G=\SO(k+1)\times \SO(n-k)\subset \SO(n+1)$ be the subgroup of block diagonal matrices that act on $\R^{n+1}=\R^{k+1}\oplus\R^{n-k}$ preserving this orthogonal splitting. Clearly, $\G$ acts on the unit sphere $\Ss^n(1)\subset\R^{n+1}$, and the unit speed geodesic $\gamma\colon \left[0,\frac\pi2\right]\to\Ss^n(1)$, given by $\gamma(\rho)=\cos \rho\, e_1+\sin \rho\, e_{n+1}$, where $\{e_j\}$ is the canonical basis of $\R^{n+1}$,  meets all $\G$-orbits in $\Ss^n(1)$ orthogonally. The orbits $\G(\gamma(0))\cong \Ss^k(1)\times\{0\}$ and $\G(\gamma(\frac\pi2))\cong \{0\}\times \Ss^{n-k-1}(1)$ are singular orbits; all other orbits $\G(\gamma(\rho))\cong \Ss^k(\cos \rho)\times \Ss^{n-k-1}(\sin \rho)$, $0<\rho<\frac\pi2$, are principal orbits. 
Using this framework, we may define $\G$-invariant metrics on $\Ss^{n}$ by specifying their values on the (open and dense) subset of principal points as the doubly warped product 
\begin{equation}\label{eq:cohom1g}\phantom{, \qquad 0<\rho<\tfrac\pi2,}
\g=\dd \rho^2+\varphi(\rho)^2\; \g_{\Ss^k} +\psi(\rho)^2 \; \g_{\Ss^{n-k-1}}, \qquad 0<\rho<\tfrac\pi2,
\end{equation}
where $\varphi$ and $\psi$ are positive functions satisfying appropriate smoothness conditions at $\rho=0$ and $\rho=\frac\pi2$, and $\g_{\Ss^d}$ is the unit round metric on $\Ss^d$. Cohomogeneity one metrics of the form \eqref{eq:cohom1g} are called \emph{diagonal}. For example, the unit round metric $\g_0=\g_{\Ss^{n}}$ is of the above form, with functions $\varphi_0(\rho)=\cos \rho$ and $\psi_0(\rho)=\sin \rho$.

The \emph{Cheeger deformation} of $\g_0$ is the $1$-parameter family $\g_s$, $s\geq0$, of diagonal cohomogeneity one metrics \eqref{eq:cohom1g} determined by the functions
\begin{equation}\label{eq:Cheeger-metric}
\varphi_s(\rho)=\frac{\cos \rho}{\sqrt{1+s\cos^2 \rho}} \quad \text{and}\quad \psi_s(\rho)=\frac{\sin \rho}{\sqrt{1+s\sin^2 \rho}},
\end{equation}
see \cite[Ex 6.46]{ABbook}. For all $s\geq0$, the metric $\g_s$ is $C^\infty$ smooth and $\G$-invariant, the orbit space of the $\G$-action on $(\Ss^n,\g_s)$ is $\Ss^n/\G=\left[0,\frac\pi2\right]$, and $\gamma$ remains a unit speed geodesic orthogonal to all $\G$-orbits. 
As the projection $\Ss^n\to\Ss^n/\G$ is distance-nonincreasing, we have
\begin{equation}\label{eq:diam-boundbelow}\phantom{, \qquad \text{ for all }s\geq0}
\diam(\Ss^n,\g_s)\geq\tfrac\pi2, \qquad \text{ for all }s\geq0.
\end{equation}
Moreover, $(\Ss^n,\g_s)$ has $\sec\geq0$ for all $s\geq0$, and it converges in Gromov--Hausdorff sense to $\Ss^n/\G=\left[0,\frac\pi2\right]$ as $s\nearrow+\infty$. 

The $\G$-orbits in $(\Ss^n,\g)$, where $\g$ is the cohomogeneity one diagonal metric \eqref{eq:cohom1g}, are isometric to the product $\G(\gamma(\rho)) = \Ss^k(\varphi(\rho))\times\Ss^{n-k-1}(\psi(\rho))$ of round spheres of radii $\varphi(\rho)$ and $\psi(\rho)$. Thus, the distance between any $x,y\in \G(\gamma(\rho))$ is
\begin{multline*}
    \dist_\g(x,y)\leq  \diam(\G(\gamma(\rho)),\g)=\sqrt{\diam(\Ss^k(\varphi(\rho)))^2+\diam(\Ss^{n-k-1}(\psi(\rho)))^2 } \\
    =\pi\sqrt{\varphi(\rho)^2+\psi(\rho)^2}.
\end{multline*}
Setting $\varphi$ and $\psi$ to be the functions in \eqref{eq:Cheeger-metric}, one easily checks that the maximum value of the above is achieved at $\rho=\frac\pi4$ for all $s\geq0$, and is equal to $\frac{\pi}{\sqrt{1+\frac{s}{2}}}$. 

The antipodal map $I\colon \Ss^n\to\Ss^n$, which acts as $I=-\Id\in\Ot(n+1)$, commutes with the $\G$-action on $(\Ss^n,\g_s)$, thus $I$ leaves invariant all $\G$-orbits. In fact, $I$ restricts to the antipodal map on each sphere factor in $\G(\gamma(\rho))$, $\rho\in \left[0,\tfrac\pi2\right]$. Thus, the displacement of $I$ on $(\Ss^n,\g_s)$ satisfies
\begin{equation*}
    \dist_{\g_s}(x,I(x))\leq\max_{\rho\in\left[0,\tfrac\pi2\right]} \diam(\G(\gamma(\rho)),\g_s)=\tfrac{\pi}{\sqrt{1+\frac{s}{2}}}.
\end{equation*}
Together with \eqref{eq:diam-boundbelow}, this proves the Theorem in the Introduction.\qed

\begin{remark}
    Not all $\G$-invariant metrics on $\Ss^n$ are diagonal, i.e., of the form \eqref{eq:cohom1g}, if $n$ is odd. For instance, let $n=3$ and $k=1$. For all $t\neq 1$, the isometry group of the Berger sphere $(\Ss^3,\g(t))$ is $\Ut(2)\subset\SO(4)$, which contains $\G=\SO(2)\times \SO(2)$, so $\g(t)$ is $\G$-invariant. However, $\g(t)$ is not of the form \eqref{eq:cohom1g} if $t\neq 1$.  Indeed, principal $\G$-orbits in $(\Ss^3,\g(t))$ are isometric to flat $2$-tori $(\G(\gamma(\rho)),\g(t))\cong \R^2/\Gamma_{(\rho,t)}$ and none of the lattices $\Gamma_{(\rho,t)}$ are rectangular if $t\neq 1$. Meanwhile, principal $\G$-orbits in $(\Ss^3,\g)$, with $\g$ as in \eqref{eq:cohom1g}, are rectangular flat tori $(\G(\gamma(\rho)),\g)\cong \R^2/2\pi\varphi(\rho)\Z\oplus2\pi\psi(\rho)\Z$.
\end{remark}

\section{Final remarks}

\subsection{Berger spheres}
\begin{figure}[!ht]
    \includegraphics[width=0.78\textwidth]{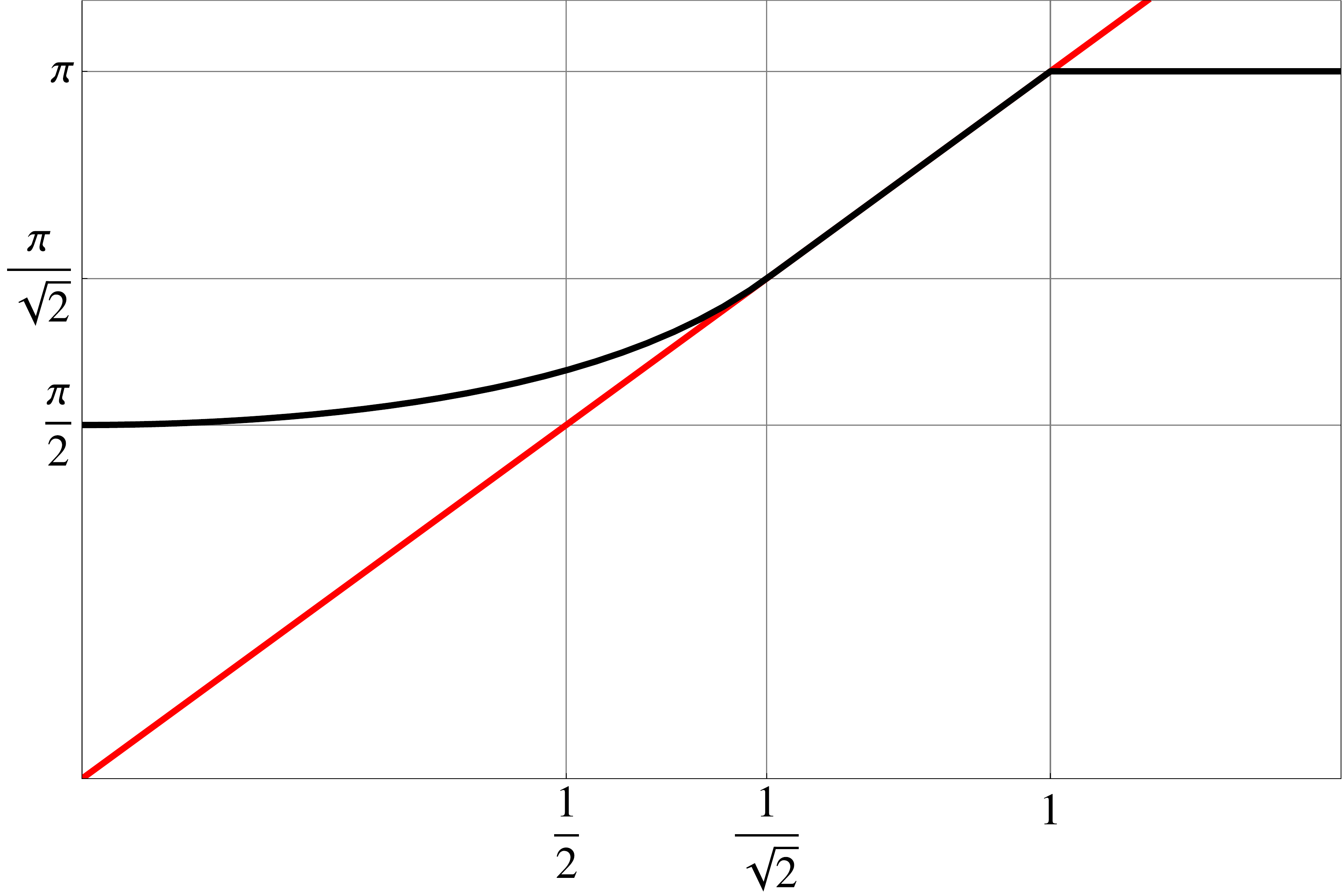}    
    \vspace{-.25cm}
    \caption{Diameter (black) and half length of Hopf circle (red) in $(\Ss^{2q+1},\g(t))$.}\label{graph}
\end{figure}

Let us expand on our discussion of the spheres $(\Ss^{2q+1},\g(t))$, whose Hopf circles are closed geodesics of length $2\pi\, t$. According to \cite{Ra85,Podobryaev17},
\begin{equation}\label{eq:diam}
\diam(\Ss^{2q+1},\g(t))=\begin{cases}
    \dfrac{\pi}{2\sqrt{1-t^2}}, & \text{ if } 0< t\leq \frac{1}{\sqrt2},\\[10pt]
    \pi \, t, & \text{ if } \frac{1}{\sqrt2}<t\leq 1,\\[6pt]
    \pi, & \text{ if }  1<t.
\end{cases}    
\end{equation}
As pairs of antipodal points $x$ and $I(x)$ are joined by half of the Hopf circle to which they belong, $\dist_{\g(t)}(x,I(x))\leq \pi\, t < \diam(\Ss^{2q+1},\g(t))$ for all $t<\frac{1}{\sqrt2}$, see \Cref{graph}.

A similar situation occurs on the Berger spheres $(\Ss^{4q+3},\mathrm h(t))$ and $(\Ss^{15},\mathrm k(t))$ obtained by scaling the unit round sphere by $t>0$ in the vertical direction of the Hopf bundles $\Ss^3\to \Ss^{4q+3}\to\mathds H P^q$ and $\Ss^7\to \Ss^{15}\to\Ss^8\!\left(\tfrac12\right)$, respectively. Namely, for all $t>0$, the projection map of these bundles remains a Riemannian submersion, and thus $\diam(\Ss^{4q+3},\mathrm h(t))\geq\diam(\mathds H P^q)=\frac\pi2$ and $\diam(\Ss^{15},\mathrm k(t))\geq\diam(\Ss^8\!\left(\tfrac12\right))=\frac\pi2$. Pairs of antipodal points belong to the same Hopf circle, hence to the same (totally geodesic) fiber, which is isometric to $\Ss^3(t)$ or $\Ss^7(t)$, so $\dist_{\g(t)}(x,I(x))\leq \pi\, t$. Thus, for sufficiently small $t>0$, these spheres also provide a negative answer to \Cref{question}.

\subsection{First Laplace eigenvalue}
Spectral geometry provides an alternative path to show that Berger spheres yield a negative answer to \Cref{question}, by considering 
\begin{equation*}
\g\longmapsto \lambda_1(M,\g)\, \diam(M,\g)^2,
\end{equation*}
where $\lambda_1(M,\g)$ is the smallest positive eigenvalue of the Laplace--Beltrami operator. This scale-invariant functional is bounded from below by $\frac{\pi^2}{4}$ on compact connected homogeneous spaces~\cite{Li80}. Moreover, $\lambda_1(\Ss^{2q+1},\g(t))\leq 4(q+1)$ for all $t>0$, since
\begin{equation*}
\lambda_1(\Ss^{2q+1}, \g(t)) = \min\big\{4(q+1),\, 2q+\tfrac{1}{t^2} \big\}
= \begin{cases}
4(q+1), &\text{if }t\leq \frac{1}{\sqrt{2q+4}},\\[2mm]
2q+\frac{1}{t^2}, &\text{if }t\geq\frac{1}{\sqrt{2q+4}},
\end{cases}
\end{equation*}
see \cite[Prop.~5.3]{BettiolPiccione13a}. Similar upper bounds on $\lambda_1$ for $(\Ss^{4q+3},\mathrm h(t))$ and $(\Ss^{15},\mathrm k(t))$ can be obtained from~\cite{BLPhomospheres}. This yields a positive diameter lower bound, independent of $t>0$, that could be used in lieu of the exact value \eqref{eq:diam} for $(\Ss^{2q+1},\g(t))$ or of the submersion lower bound $\frac\pi2$ in general. However, this spectral lower bound on the diameter is weaker than the latter, and becomes arbitrarily small as $q\nearrow+\infty$.

\bibliographystyle{plain}

\begin{thebibliography}{BLP22}

\bibitem[AB15]{ABbook}
{\sc M.~M. Alexandrino and R.~G. Bettiol}.
\newblock {\em Lie groups and geometric aspects of isometric actions}.
\newblock Springer, Cham, 2015.
DOI: \href{http://dx.doi.org/10.1007/978-3-319-16613-1} {10.1007/978-3-319-16613-1}.

\bibitem[BP13]{BettiolPiccione13a}
	{\sc R.~G.~Bettiol, P.~Piccione}.
	{\it Bifurcation and local rigidity of homogeneous solutions to the {Y}amabe problem on spheres.}
	Calc. Var. Partial Differential Equations \textbf{47}:3--4 (2013), 789--807.
	DOI: \href{http://dx.doi.org/10.1007/s00526-012-0535-y} {10.1007/s00526-012-0535-y}.


\bibitem[BLP22]{BLPhomospheres}
	{\sc R.~G.~Bettiol, E.~A.~Lauret, P.~Piccione}.
	{\it The first eigenvalue of a homogeneous CROSS.}
	J. Geom. Anal. \textbf{32} (2022), 76.
	DOI: \href{https://doi.org/10.1007/s12220-021-00826-7} {10.1007/s12220-021-00826-7}.

\bibitem[BH99]{bh}
    {\sc M.~Bridson, A.~Haefliger}.
    {\it Metric spaces of non-positive curvature}.
    Grundlehren Math. Wiss., 319. Springer-Verlag (1999).
    DOI: \href{http://dx.doi.org/10.1007/978-3-662-12494-9} {10.1007/978-3-662-12494-9}.


\bibitem[GP93]{GrovePetersen}
	{\sc K.~Grove, P.~Petersen}.
	{\it A radius sphere theorem}.
	Invent. Math. \textbf{112}:3 (1993), 577--583.
	DOI: \href{http://dx.doi.org/10.1007/BF01232447} {10.1007/BF01232447}.


\bibitem[Li80]{Li80}
	{\sc P.~Li}.
	{\it Eigenvalue estimates on homogeneous manifolds.}
	Comment.\ Math.\ Helvetici \textbf{55} (1980), 347--363.
	DOI: \href{http://dx.doi.org/10.1007/BF02566692} {10.1007/BF02566692}.

\bibitem[Ni01]{Nikonorov01diameter}
	{\sc Yu.~G.~{Nikonorov}}.
	{\it On the geodesic diameter of surfaces possessing an involutory isometry} (Russian).
	Tr. Rubtsovsk. Ind. Inst. \textbf{9} (2001), 62--65. 
	English translation:
	\href{http://arxiv.org/abs/1811.01173} {arXiv:1811.01173}. 
	

\bibitem[Po18a]{Podobryaev17}
	{\sc A.~V.~Podobryaev}.
	{\it Diameter of the Berger Sphere}.
	Math. Notes \textbf{103}:5--6 (2018), 846--851. 
    DOI: \href{http://dx.doi.org/10.1134/S0001434618050188} {10.1134/S0001434618050188}.

\bibitem[Po18b]{Podobryaev18}
	{\sc A.~V.~Podobryaev}.
	{\it Antipodal points and diameter of a sphere}.
	Russ. J. Nonlinear Dyn. \textbf{14}:4 (2018), 579--581. 
    DOI: \href{http://dx.doi.org/10.20537/nd180410} {10.20537/nd180410}.
    
\bibitem[Ra85]{Ra85}
    {\sc C.~Rakotoniaina}.
    {\it Cut locus of the B-spheres}.
    Ann. Global Anal. Geom. \textbf{3}:3 (1985), 313--327. 
    DOI: \href{http://dx.doi.org/10.1007/BF00130483} {10.1007/BF00130483}.


\end{thebibliography}

\end{document}